\documentclass[11pt]{article}
\usepackage{amsfonts, amsmath, amssymb, amsthm}

\PassOptionsToPackage{obeyspaces}{url}
\usepackage[colorlinks=true,citecolor=blue,urlcolor=blue,linkcolor=blue,bookmarksopen=true]{hyperref}
\usepackage{breakurl}
\usepackage{tikz}
\usetikzlibrary{arrows}
\usetikzlibrary{calc}
\usepackage{enumitem}

\usepackage{etoolbox} 
\patchcmd{\thebibliography}{\leftmargin\labelwidth}{\leftmargin\labelwidth\addtolength\itemsep{-0.208\baselineskip}}{}{}

\author{Boris Bukh\thanks{Department of Mathematical Sciences, Carnegie Mellon University, Pittsburgh, PA, USA. \texttt{bbukh@math.cmu.edu}. Supported in part by Sloan Research Fellowship and by U.S.\ taxpayers through NSF CAREER grant DMS-1555149.}\and 
Christopher Cox\thanks{Department of Mathematical Sciences, Carnegie Mellon University, Pittsburgh, PA, USA. \texttt{cocox@andrew.cmu.edu}. Supported in part by U.S.\ taxpayers through NSF CAREER grant DMS-1555149.}}

\title{Nearly orthogonal vectors and small antipodal spherical codes}
\date{\today}

\newtheorem{theorem}{Theorem}
\newtheorem{lemma}[theorem]{Lemma}
\newtheorem{cor}[theorem]{Corollary}
\newtheorem{prop}[theorem]{Proposition}
\newtheorem{conj}[theorem]{Conjecture}
\theoremstyle{definition}
\newtheorem{defn}[theorem]{Definition}
\newtheorem{fact}{Fact}

\newcommand*{\abs}[1]{\lvert #1\rvert}                           

\newcommand{\C}{\mathbb{C}}
\newcommand{\Z}{\mathbb{Z}}

\newcommand{\E}{\mathbb{E}}
\newcommand{\R}{\mathbb{R}}
\renewcommand{\P}{\mathcal{P}}
\newcommand{\bh}{\mathbb{H}}

\newcommand{\inner}{\theta}
\newcommand{\tr}{\operatorname{tr}}
\newcommand{\rank}{\operatorname{rk}}
\newcommand{\supp}{\operatorname{supp}}
\newcommand{\Lone}{\mathcal{L}}
\newcommand{\SLone}{\mathcal{SL}}
\newcommand{\off}{\operatorname{off}}
\newcommand{\Span}{\operatorname{span}}
\DeclareMathOperator{\GL}{GL}

\begin{document}

	\maketitle
\begin{abstract}
	How can $d+k$ vectors in $\R^d$ be arranged so that they are as close to orthogonal as possible? In particular, define $\inner(d,k):=\min_X\max_{x\neq y\in X}|\langle x,y\rangle|$ where the minimum is taken over all collections of $d+k$ unit vectors $X\subseteq\R^d$. In this paper, we focus on the case where $k$ is fixed and $d\to\infty$. In establishing bounds on $\theta(d,k)$, we find an intimate connection to the existence of systems of ${k+1\choose 2}$ equiangular lines in $\R^k$. Using this connection, we are able to pin down $\theta(d,k)$ whenever $k\in\{1,2,3,7,23\}$ and establish asymptotics for general $k$. The main tool is an upper bound on $\E_{x,y\sim\mu}|\langle x,y\rangle|$ whenever $\mu$ is an isotropic probability mass on $\R^k$, which may be of independent interest. Our results translate naturally to the analogous question in $\C^d$. In this case, the question relates to the existence of systems of $k^2$ equiangular lines in $\C^k$, also known as SIC-POVM in physics literature.
\end{abstract}

\section{Introduction}
How can a given number of points be arranged on a sphere in $\R^d$ so that they are as far from each other as possible? This is
a basic problem in coding theory; for example, the book \cite{codes_book} is devoted to this problem exclusively. Such point arrangements
are called \emph{spherical codes}. Most constructions of spherical codes are symmetric. 
Here we consider the \emph{antipodal codes}, in which the points come in pairs $x,-x$.
In other words, we seek arrangements of $d+k$ unit vectors in $\R^d$ so that they are as close to orthogonal as possible. An alternative point of view is that these are codes in the projective space $\mathbb{RP}^{d-1}$; for example, see \cite{CohnKumarProjective}. We focus on the case when $k$ is small.

As we will see, this question relates to the problem of the existence of large families of equiangular lines in~$\R^k$.
Similarly, the analogous question for unit vectors in $\C^d$ relates to equiangular lines in $\C^k$, which
are the mathematical underpinning of symmetric informationally complete measurements in quantum theory~\cite{SIC}. Because of this,
we elect to treat the real and complex cases in parallel. Henceforth, we denote by $\bh$ the underlying field, which 
can be either~$\R$ or~$\C$.

For $\bh\in\{\R,\C\}$, define the parameter
\[
\inner_\bh(d,k):=\min_X\max_{x\neq y\in X}|\langle x,y\rangle|,
\]
where the minimum is taken over all collections of $d+k$ unit vectors $X\subseteq\bh^d$. In this paper, we prove bounds on $\inner_\bh(d,k)$ when $k$ is fixed and $d\to\infty$.

For a collection of vectors $X=\{x_1,\dots,x_n\}\subseteq\bh^d$, the \emph{Gram matrix} is the matrix $A\in\bh^{n\times n}$ where $A_{ij}=\langle x_i,x_j\rangle$. It will be easier to work with the Gram matrices than with the vectors themselves. 

For a matrix $A\in\bh^{n\times n}$, define $\off(A):=\max_{i\neq j}|A_{ij}|$. By considering Gram matrices, one can equivalently define $\inner_\bh(d,k)=\min_A\off(A)$ where the minimum is taken over all $A\in\bh^{(d+k)\times(d+k)}$ with $\rank(A)=d$ where $A_{ii}=1$ for every $i$ and $A$ is Hermitian and positive semidefinite. Our techniques are not specialized to Hermitian, positive semidefinite matrices, so we also define 
\[
\off_\bh(d,k):=\min_A\off(A),
\]
where the minimum is taken over all $A\in\bh^{(d+k)\times(d+k)}$ with $\rank(A)=d$ and $A_{ii}=1$ for every $i$. Note that $\off_\bh(d,k)\leq\inner_\bh(d,k)$.

In Section~\ref{sec:lower}, we establish lower bounds on $\off_\bh(d,k)$, and in Section~\ref{sec:upper}, we give constructions to yield upper bounds on $\inner_\bh(d,k)$. Throughout both of these sections, we will show an intimate connection between determining these parameters and the existence of large systems of equiangular lines in $\bh^k$.

\begin{defn}
	A system of equiangular lines is a collection of lines through the origin which pairwise meet at the same angle.
	We identify a line with a spanning unit vector, so formally, a system of equiangular lines in $\bh^k$ is a collection of unit vectors $X\subseteq\bh^k$ so that there is some $\beta\in\R$ where $|\langle x,y\rangle|=\beta$ for all $x\neq y\in X$. 
	
	It is known that if $X\subseteq\R^k$ is a system of equiangular lines, then $|X|\leq{k+1\choose 2}$ and if $X\subseteq\C^k$ is a system of equiangular lines, then $|X|\leq k^2$.
\end{defn}

The main results of this paper are as follows:
\begin{theorem}\label{thm:genlower}\hspace{2em}\vspace{-1ex}
	\begin{enumerate}[label=(\alph*)]
		\item For positive integers $d,k$,
		\[
		\off_\R(d,k)\geq{1\over \alpha_k(d+k)-1},
		\]
		where $\alpha_k={(k-1)\sqrt{k+2}+2\over k(k+1)}$. If equality holds, then there exists a system of ${k+1\choose 2}$ equiangular lines over $\R^k$ and $d\equiv -k\pmod{{k+1\choose 2}}$.
		
		\item For positive integers $d,k$,
		\[
		\off_\C(d,k)\geq{1\over\alpha_k^*(d+k)-1},
		\]
		where $\alpha_k^*={(k-1)\sqrt{k+1}+1\over k^2}$. If equality holds, then there exists a system of $k^2$ equiangular lines over $\C^k$ and $d\equiv -k\pmod{k^2}$.
	\end{enumerate}
\end{theorem}
This is an improvement over the classical Welch bound (which is recalled as Theorem~\ref{thm:firstlower} below) when $k\leq O(d^{1/2})$. 
It is a quantitative improvement of a result of Cohn--Kumar--Minton \cite[Corollary~2.13]{CohnKumarProjective}
which asserts that Welch bound is not sharp for $k\leq O(d^{1/2})$, without providing a better bound.

A computer-assisted proof of the case $(d,k)=(4,2)$ of Theorem~\ref{thm:genlower} was recently given by Fickus--Jasper--Mixon \cite{fickus_jasper_mixon}.

The above theorem will follow as a corollary of Theorems~\ref{thm:probbound} and~\ref{thm:isobound}, which will be proved in Section~\ref{sec:lower}. Furthermore, the following theorem, which will be proved in Section~\ref{sec:upper}, will show that equality does, in fact, hold under the stated conditions.

\begin{theorem}\label{thm:equiconst}\hspace{2em}\vspace{-1ex}
	\begin{enumerate}[label=(\alph*)]
		\item If there is a system of ${k+1\choose 2}$ equiangular lines in $\R^k$ and $d\equiv -k\pmod{{k+1\choose 2}}$, then 
		\[
		\off_\R(d,k)=\inner_\R(d,k)={1\over\alpha_k(d+k)-1},
		\]
		where $\alpha_k={(k-1)\sqrt{k+2}+2\over k(k+1)}$.
		
		\item If there is a system of $k^2$ equiangular lines in $\C^k$ and $d\equiv -k\pmod {k^2}$, then
		\[
		\off_\C(d,k)=\inner_\C(d,k)={1\over\alpha_k^*(d+k)-1},
		\]
		where $\alpha_k^*={(k-1)\sqrt{k+1}+1\over k^2}$.
	\end{enumerate}
\end{theorem}

The usual way of proving bounds on codes is to use linear programming. In the context of spherical codes, the relevant linear
program first appeared in the work of Delsarte, Goethals and Seidel \cite{Delsarte1977}. See \cite[Chapter~2]{codes_book} for the general exposition,
and \cite{astola} for the case of few vectors.

In contrast, we establish Theorem~\ref{thm:genlower} by relating the problem to that of bounding the first moment of isotropic measures.

For a probability mass $\mu$, we write $\E_{x\sim\mu}f(x)$ to denote the expected value of $f(x)$ where $x$ is distributed according to $\mu$.
We also use $\E_{x,y\sim\mu}f(x,y):=\E_{x\sim\mu}\E_{y\sim\mu}f(x,y)$.
\begin{defn}
	For $\bh\in\{\R,\C\}$, a probability mass $\mu$ on $\bh^k$ is called \emph{isotropic} if $\E_{x\sim\mu}|\langle x,v\rangle|^2={1\over k}\Vert v\Vert^2$ for every $v\in\bh^k$. Equivalently, $\mu$ is isotropic if $\E_{x\sim\mu}xx^*={1\over k}I_k$. Such a probability mass is also called a \emph{probabilistic tight frame} with frame constant $1/k$ (see~\cite{Ehler2013} for a survey).
\end{defn}

We show the following:
\begin{lemma}\label{lem:iso}\hspace{2em}
	\begin{enumerate}[label=(\alph*)]
	\item If $\mu$ is an isotropic probability mass on $\R^k$, then 
	\[
	\E_{x,y\sim\mu}|\langle x,y\rangle|\leq {(k-1)\sqrt{k+2}+2\over k(k+1)},
	\]
	with equality if and only if there exists $X\subseteq\R^k$, a system of ${k+1\choose 2}$ equiangular lines, and $\mu$ satisfies $\mu(x)+\mu(-x)=1/{k+1\choose 2}$ for every $x\in X$.

	\item If $\mu$ is an isotropic probability mass on $\C^k$, then 
	\[
	\E_{x,y\sim\mu}|\langle x,y\rangle|\leq{(k-1)\sqrt{k+1}+1\over k^2},
	\]
	with equality if and only if there exists $X\subseteq\C^k$, a system of $k^2$ equiangular lines, and $\mu$ satisfies $\mu(x)+\mu(-x)=1/k^2$ for every $x\in X$.
	\end{enumerate}
	
\end{lemma}

Theorem~\ref{thm:isobound} shows the connection between the above lemma and Theorem~\ref{thm:genlower}.

As there are systems of ${k+1\choose 2}$ equiangular lines over $\R^k$ whenever $k\in\{1,2,3,7,23\}$, we can give tight answers for infinitely many $d$ in these cases; see Corollary~\ref{cor:equicases} for the exact values. See \cite{Glazyrin2016,greaves_koolen_munemasa_szollosi,yu_new_bounds} for the known bounds of the size of the largest system of equiangular lines in $\R^k$.

Even in the cases not covered by Theorem~\ref{thm:equiconst}, we will still show that Theorem~\ref{thm:genlower} is asymptotically tight.
\begin{theorem}\label{thm:genupper}
	Let $\bh\in\{\R,\C\}$. For every $\epsilon>0$, there is an integer $k_0$ so that for any fixed $k\geq k_0$,
	\[
	\inner_\bh(d,k)\leq\bigl(1+o(1)\bigr){(1+\epsilon)\sqrt{k}\over d},
	\]
	where $o(1)\to 0$ as $d\to\infty$.
\end{theorem}

The above theorem will be established in multiple parts. First, Theorem~\ref{thm:offby2} will show that $\inner_\R(d,k)\leq(1+o(1)){\sqrt{k+4}\over d}$ whenever $k$ is a power of $4$ and show that $\inner_\R(d,k)\leq(1+o(1)){2\sqrt{k+1}\over d}$ for general $k$. Theorem~\ref{thm:asymtightcomp} will establish Theorem~\ref{thm:genupper} in the case of complex numbers and show that in this case we can take $k_0=O(\epsilon^{-40/19})$. Finally, Theorem~\ref{thm:genupper} will be established fully in the case of the reals by Theorem~\ref{thm:realasym}.

\paragraph{Acknowledgments.} We thank William Martin for inspiring discussions. We also thank
organizers of the Ninth Discrete Geometry and Algebraic Combinatorics Conference, where these discussions took place. The conference 
was supported by NSF grant DMS-162360. We also thank the referee for helpful suggestions.

\section{Lower bounds}\label{sec:lower}

\paragraph{Basic bound and the case $k=1$.}
We begin with a simple lower bound which is originally due to Welch~\cite{welch1974lower} and has been rediscovered various times in the literature, for example \cite[Lemma~2.2]{Alon2009} and \cite[Lemma~3.2]{YNLunpub}. We give a proof for completeness.

\begin{theorem}\label{thm:firstlower}
	For $\bh\in\{\R,\C\}$, if $d,k$ are positive integers, then $\off_\bh(d,k)\geq\sqrt{{k\over d(d+k-1)}}$.
\end{theorem}
\begin{proof}
	Let $A\in\bh^{(d+k)\times(d+k)}$ with $1$'s on the diagonal and $\rank(A)\leq d$. Then
	\[
	\tr (A^*A)=\sum_{i,j}|A_{ij}|^2=(d+k)+\sum_{i\neq j}|A_{ij}|^2\leq (d+k)+(d+k)(d+k-1)\off(A)^2.
	\]
	On the other hand, $\tr (A^*A)\geq |\tr( A)|^2/\rank(A)$ (see Proposition~\ref{prop:infinitematrix} for a proof), so
	\[
	(d+k)+(d+k)(d+k-1)\off(A)^2\geq\tr (A^*A)\geq{|\tr (A)|^2\over d}={(d+k)^2\over d}.
	\]
	Rearranging these inequalities yields $\off(A)\geq\sqrt{{k\over d(d+k-1)}}$, so the same bound holds for $\off_\bh(d,k)$.
\end{proof}

Before moving on, we note that the above observation suffices to determine $\off_\bh(d,1)$ and $\inner_\bh(d,1)$.

\begin{cor}\label{cor:k=1}
	For $\bh\in\{\R,\C\}$ and for any positive integer $d$, $\off_\bh(d,1)=\inner_\bh(d,1)={1\over d}$.
\end{cor}
\begin{proof}
	The lower bound follows from Theorem~\ref{thm:firstlower}. For the upper bound, let $x_1,\dots,x_{d+1}$ be the vertices of a unit regular simplex in $\R^d$ centered at the origin. Then for all $i\neq j$, we have $\langle x_i,x_j\rangle=-{1\over d}$, so $\theta_\R(d,1)\leq{1\over d}$. As $\theta_\C(d,1)\leq\theta_\R(d,1)$, this establishes the claim.
\end{proof}

\paragraph{Connection to isotropic measures.} We now turn our attention to the general case. Throughout the following, whenever we discus a probability mass $\mu$ on $\bh^k$, $\mu$ will be assumed to be Borel. For such a $\mu$, we use $\E_{x\sim\mu} f(x)$ to denote the expected value of the function $f$ where $x$ is distributed according to $\mu$. We also use $\E_{x,y\sim\mu} f(x,y):=\E_{x\sim \mu}\E_{y\sim\mu} f(x,y)$. When the probability mass $\mu$ is understood, we will omit writing it. Recall that the \emph{support} of $\mu$, denoted $\supp(\mu)$, is the collection of all $x\in\bh^k$ for which every ball centered at $x$ has positive mass.

The following parameter will play a crucial role in our bounds. 
\begin{defn}
	For $\bh\in\{\R,\C\}$, let $\mu$ be a nonzero probability mass on $\bh^k$ and define
	\[
	\Lone_\bh(\mu):=\inf_{y\in\supp(\mu)\setminus\{0\}}\inf_{v\in\bh^k\setminus\{0\}}{\E_{x\sim\mu}|\langle v,x\rangle|\over|\langle v,y\rangle|}.
	\]
\end{defn}

We care about the parameter $\Lone_\bh(\mu)$ only when $\mu$ is of a certain form. Define $\P_\bh(d,k)$ to be the collection of all probability masses $\mu$ on $\bh^k$ for which there is a (multi)set $X$ of $d+k$ vectors over $\bh^k$ with $\Span(X)=\bh^k$ and $\mu$ is the uniform distribution over $X$. In other words, $\P_\bh(d,k)$ is the collection of all probability masses $\mu$ where $\supp(\mu)$ is finite, $\supp(\mu)$ spans $\bh^k$ and $(d+k)\mu(x)\in\Z$ for all $x\in\supp(\mu)$.

We then define
\[
\SLone_\bh(d,k):=\sup_{\mu\in\P_\bh(d,k)}\Lone_\bh(\mu).
\]
Proposition~\ref{prop:supachieved} will show that we may replace the above supremum with a maximum.

\begin{theorem}\label{thm:probbound}
	For $\bh\in\{\R,\C\}$, if $d,k$ are positive integers, then 
	\[
	\off_\bh(d,k)\geq{1\over\SLone_\bh(d,k)(d+k)-1}.
	\]
\end{theorem}
\begin{proof}
	Let $A\in\bh^{(d+k)\times(d+k)}$ with $1$'s on the diagonal and $\rank(A)\leq d$. Thus $\dim\ker A\geq k$, so there is some $N\in\bh^{(d+k)\times k}$ with $\rank(N)=k$ and $AN=0$. Let $y_i$ be the $i$th \emph{row} of $N$, so we have $\bigl(\langle v,y_1\rangle,\langle v,y_2\rangle,\dots,\langle v,y_{d+k}\rangle\bigr)^T\in\ker A$ for every $v\in\bh^k$. Thus, for any fixed $i\in[d+k]$,
	\[
	0=\sum_j A_{ij}\langle v,y_j\rangle=\langle v,y_i\rangle+\sum_{j\neq i}A_{ij}\langle v,y_j\rangle,
	\]
	so,
	\[
	|\langle v,y_i\rangle|=\bigg|\sum_{j\neq i}A_{ij}\langle v,y_j\rangle\bigg|\leq\off(A)\sum_{j\neq i}|\langle v,y_j\rangle|.
	\]
	Solving for $\off(A)$, if $\langle v,y_i\rangle\neq 0$,
	\[
	\off(A)\geq{|\langle v,y_i\rangle|\over\sum_{j\neq i}|\langle v,y_j\rangle|}=\biggl({1\over|\langle v,y_i\rangle|}\sum_j|\langle v,y_j\rangle|-1\biggr)^{-1}.
	\]
	As this bound holds for all $i\in[d+k]$ and $v\in\bh^k$ with $\langle v,y_i\rangle\neq 0$, if $\mu$ is the uniform distribution over the (multi)set $\{y_1,\dots,y_{d+k}\}$, we have
	\[
	\off(A)\geq\sup_{y\in\supp(\mu)\setminus\{0\}}\sup_{v\in\bh^k\setminus\{0\}}\biggl({\E_x|\langle v,x\rangle|\over|\langle v,y\rangle|}(d+k)-1\biggr)^{-1}={1\over\Lone_\bh(\mu)(d+k)-1}.
	\]
	Finally, as $\{y_1,\dots,y_{d+k}\}\subseteq\bh^k$ and $\rank(N)=k$, we know that $\Span\{y_1,\dots,y_{d+k}\}=\bh^k$, and so $\mu\in\P_\bh(d,k)$. As such, $\Lone_\bh(\mu)\leq\SLone_\bh(d,k)$, implying
	\[
	\off(A)\geq{1\over\SLone_\bh(d,k)(d+k)-1},
	\]
	which yields the same lower bound on $\off_\bh(d,k)$.
\end{proof}

Thus, in order to obtain lower bounds on $\off_\bh(d,k)$, it suffices to establish upper bounds on~$\SLone_\bh(d,k)$.

For a matrix $Q\in\GL_k(\bh)$ and a probability mass $\mu$ on $\bh^k$, let $Q\mu$ be the probability mass defined by $Q\mu (S):=\mu (Q^{-1}S)$ for every Borel set $S$. Recalling that
$\mu$ is isotropic if $\E_{x\sim \mu} xx^*=\frac{1}{k}I_k$, it is not difficult to see that if $\mu$ is a probability mass on $\bh^k$, then $\supp(\mu)$ spans $\bh^k$ if and only if there is some $Q\in\GL_k(\bh)$ for which $Q\mu $ is isotropic.


The following proposition shows that, when considering $\Lone_\bh(\mu)$, we may always suppose that $\mu$ is isotropic.

\begin{prop}\label{prop:iso}
	If $\mu$ is a probability mass on $\bh^k$ and $Q\in\GL_k(\bh)$, then $\Lone_\bh(\mu)=\Lone_\bh(Q\mu )$. 
\end{prop}
\begin{proof}
	For any $y\in\supp(Q\mu)\setminus\{0\}$ and $v\in\bh^k\setminus\{0\}$, we find
	\[
	{\E_{x\sim Q\mu }|\langle x,v\rangle|\over|\langle y,v\rangle|}={\E_{x\sim\mu} |\langle Qx,v\rangle|\over|\langle QQ^{-1}y,v\rangle|}={\E_{x\sim\mu}|\langle x,Q^*v\rangle|\over|\langle Q^{-1}y,Q^*v\rangle|}.
	\]
	As $\supp(Q\mu)=Q\supp(\mu)$, this establishes the claim.
\end{proof}

\paragraph{First moment of isotropic measures.} We now focus on proving Lemma~\ref{lem:iso}, which will be key in establishing upper bounds on $\SLone_\bh(d,k)$. To do so, we will need two facts about ``infinite matrices''.

Let $\Omega$ be a set and $f\colon\Omega^2\to\bh$. The \emph{rank} of $f$, denoted $\rank (f)$, is defined to be the smallest $r$ for which there are functions $g_i,h_i\colon\Omega\to\bh$, $i\in[r]$, so that $f(x,y)=\sum_{i=1}^r g_i(x)h_i(y)$ for every $x,y\in \Omega$. If there is no such $r$, define $\rank(f)=\infty$. Notice that if $|\Omega|<\infty$, then the rank of $f$ is the rank of the matrix $A$ defined by $A_{xy}=f(x,y)$.
Let $f^*$ be defined by $f^*(x,y)=\overline{f(y,x)}$ and $\overline{f}$ be defined by $\overline{f}(x,y)=\overline{f(x,y)}$. The following inequality will be essential in the proof of Lemma~\ref{lem:iso}.

\begin{prop}\label{prop:infinitematrix}
	For $\bh\in\{\R,\C\}$, let $f\colon\Omega^2\to\bh$ and $\mu$ be a probability mass on $\Omega$. If $\rank(f)<\infty$, then
	\[
	\E_{x,y\sim\mu}f^*(x,y)f(x,y)\geq{|\E_{x\sim\mu}f(x,x)|^2\over \rank(f)}.
	\]
\end{prop}
\begin{proof}
	For completeness, we first give a proof when $|\Omega|<\infty$ and $\mu$ is the uniform distribution over $\Omega$. In this case, let $A$ be the matrix with $A_{x,y}=f(x,y)$. Let $\lambda_1,\dots,\lambda_{\rank(A)}$ be the nonzero eigenvalues of $A$ and $\sigma_1,\dots,\sigma_{\rank(A)}$ be the nonzero singular values of $A$. It is well-known that $\sum_{i=1}^{\rank(A)}|\lambda_i|\leq\sum_{i=1}^{\rank(A)}\sigma_i$ (see \cite[Eq.~(II.23)]{bhatia}). Therefore, by Cauchy--Schwarz,
	\[
	\tr(A^*A)=\sum_{i=1}^{\rank(A)}\sigma_i^2\geq{1\over\rank(A)}\biggl(\sum_{i=1}^{\rank(A)}\sigma_i\biggr)^2\geq{1\over\rank(A)}\biggl(\sum_{i=1}^{\rank(A)}|\lambda_i|\biggr)^2\geq{|\tr(A)|^2\over\rank(A)}.
	\]
	
	Now, for a general $\Omega$ and $\mu$, let $x_1,\dots,x_n$ be independent samples from $\Omega$ according to $\mu$. If $f'$ denotes the restriction of $f$ to $\{x_1,\dots,x_n\}^2$, then certainly $\rank(f')\leq\rank(f)$. Hence, from above,
	\[
	{1\over n^2}\sum_{i,j}f^*(x_i,x_j)f(x_i,x_j)\geq {1\over \rank(f)}\biggl\lvert{1\over n}\sum_if(x_i,x_i)\biggr\rvert^2.
	\]
        Taking the expectation of both sides over the random choice of the samples $x_1,\dotsc,x_n$, and using that $\E[X^2]\geq \E[X]^2$ for any random variable $X$, we obtain
        \[
          \frac{n-1}{n}\E_{x,y\sim\mu}f^*(x,y)f(x,y)+\frac{1}{n}\E_{x\sim \mu} \abs{f(x,x)}^2\geq \frac{1}{\rank(f)}\bigl\lvert \E_{x\sim \mu} f(x,x)\bigr\rvert^2. 
        \]
        Taking the limit $n\to\infty$ establishes the claim.
\end{proof}

We will require also the following observation, which generalizes the corresponding property of Hadamard products.

\begin{prop}\label{prop:hadamard}
	For $\bh\in\{\R,\C\}$, let $f\colon\Omega^2\to\bh$. If $\rank(f)=r$, then $\rank(f^2)\leq{r+1\choose 2}$ and $\rank(\overline{f}f)\leq r^2$.
\end{prop}
\begin{proof}
	Let $g_i,h_i\colon\Omega\to\bh$, $i\in[r]$, be such that $f(x,y)=\sum_{i=1}^r g_i(x)h_i(y)$ for every $x,y\in\Omega$. As such, 
	\[
	f(x,y)^2=\sum_{i,j}g_i(x)g_j(x)h_i(y)h_j(y)=\sum_{i\leq j}g'_{ij}(x)h'_{ij}(y),
	\]
	where $g'_{ii}=g_i^2$, $h'_{ii}=h_i^2$, and for $i<j$, $g'_{ij}={1\over 2}g_ig_j$ and $h'_{ij}={1\over 2}h_ih_j$. Therefore, $\rank(f^2)\leq{r+1\choose 2}$.
	
	Similarly,
	\[
	\overline{f}(x,y)f(x,y)=\sum_{i,j}\overline{g_i}(x)g_j(x)\overline{h_i}(y)h_j(y)=\sum_{i,j}g'_{ij}(x)h'_{ij}(y),
	\]
	where $g'_{ij}=\overline{g_i}g_j$ and $h'_{ij}=\overline{h_i}h_j$, so $\rank(\overline{f}f)\leq r^2$.
\end{proof}

\begin{proof}[Proof of Lemma~\ref{lem:iso}]
	We first establish the upper bound. For $\bh\in\{\R,\C\}$, let $\mu$ be an isotropic probability mass on $\bh^k$. The cases where $\bh=\R$ and $\bh=\C$ will be almost identical. We will break into cases when necessary.
	
	As a technical detail, we must first assure that $\Pr_\mu[x=0]=0$. To do this, set $p=1-\Pr_\mu[x=0]$ and notice that $p>0$ as $\supp(\mu)$ spans $\bh^k$. Let $\mu'$ be the probability mass which is $\mu$ conditioned on the event $\{x\neq 0\}$. We notice that
	\[
	\E_{x\sim\mu'}|\langle x,v\rangle|^2={1\over p}\E_{x\sim\mu}|\langle x,v\rangle|^2\text{ for every $v\in\bh^k$},\quad\text{and}\quad \E_{x,y\sim\mu'}|\langle x,y\rangle|={1\over p^2}\E_{x,y\sim\mu}|\langle x,y\rangle|.
	\]
	Therefore, if $Q=\sqrt{p}I_k$, then $Q\mu'$ is isotropic and $\E_{x,y\sim Q\mu'}|\langle x,y\rangle|={1\over p}\E_{x,y\sim\mu}|\langle x,y\rangle|$. If $p<1$, then $\E_{x,y\sim Q\mu'}|\langle x,y\rangle|>\E_{x,y\sim\mu}|\langle x,y\rangle|$, and so we may replace $\mu$ by $Q\mu'$ and upper bound $\E_{x,y\sim Q\mu'}|\langle x,y\rangle|$. Hence, we may assume that
$\Pr_\mu[x=0]=0$ in what follows.
	
	From now on, we will compress notation and write $\E_x$ in lieu of $\E_{x\sim \mu}$. Set \[\alpha:=\E_{x,y}|\langle x,y\rangle|.\]
	
	For $\beta\geq 0$, we will establish upper and lower bounds on
	\[
	M(\beta):=\E_{x,y}\biggl({|\langle x,y\rangle|\over\sqrt{\Vert x\Vert\Vert y\Vert}}-\beta\sqrt{\Vert x\Vert\Vert y\Vert}\biggr)^2,
	\]
	which is well-defined since $\Pr[x=0]=0$.
	For the upper bound, we begin by expanding
	\[
	M(\beta) = \E_{x,y}{|\langle x,y\rangle|^2\over\Vert x\Vert\Vert y\Vert}-2\beta\E_{x,y}|\langle x,y\rangle|+\beta^2\E_{x,y}\Vert x\Vert\Vert y\Vert.
	\]
	By Cauchy--Schwarz, recalling that $\E_{x}|\langle x,y\rangle|^2={1\over k}\Vert y\Vert^2$ for any $y\in\bh^k$, we obtain
	\begin{equation}\label{eqn:norm}
	\E_{x,y}{|\langle x,y\rangle|^2\over\Vert x\Vert\Vert y\Vert}\leq\sqrt{\E_{x,y}{|\langle x,y\rangle|^2\over\Vert x\Vert^2}}\sqrt{\E_{x,y}{|\langle x,y\rangle|^2\over\Vert y\Vert^2}}=\sqrt{{1\over k}\E_x{\Vert x\Vert^2\over\Vert x\Vert^2}}\sqrt{{1\over k}\E_y{\Vert y\Vert^2\over\Vert y\Vert^2}}={1\over k}.
	\end{equation}
	Therefore,
	\begin{equation}\label{eqn:upper}
	M(\beta)\leq{1\over k}-2\beta\alpha+\beta^2\bigl(\E_x\Vert x\Vert\bigr)^2.
	\end{equation}

	For the lower bound, we first write,
	\begin{align}\label{eqn:CS}
		M(\beta) &= \E_{x,y}\biggl(\biggl({|\langle x,y\rangle|\over\sqrt{\Vert x\Vert\Vert y\Vert}}-\beta\sqrt{\Vert x\Vert\Vert y\Vert}\biggr){|\langle x,y\rangle|+\beta\Vert x\Vert\Vert y\Vert\over|\langle x,y\rangle|+\beta\Vert x\Vert\Vert y\Vert}\biggr)^2 \nonumber\\
				 &= \E_{x,y}\biggl({|\langle x,y\rangle|^2-\beta^2\Vert x\Vert^2\Vert y\Vert^2\over\sqrt{\Vert x\Vert\Vert y\Vert}\bigl(|\langle x,y\rangle|+\beta\Vert x\Vert\Vert y\Vert\bigr)}\biggr)^2 \nonumber\\
				 &\geq\E_{x,y}\biggl({|\langle x,y\rangle|^2-\beta^2\Vert x\Vert^2\Vert y\Vert^2\over(1+\beta)\Vert x\Vert^{3/2}\Vert y\Vert^{3/2}}\biggr)^2,
	\end{align}
	where the last inequality follows by applying Cauchy--Schwarz to the denominator.

	Set $\Omega=\bh^k\setminus\{0\}$, and define $f\colon\Omega^2\to\bh$ by 
	\[
	f(x,y):={|\langle x,y\rangle|^2-\beta^2\Vert x\Vert^2\Vert y\Vert^2\over(1+\beta)\Vert x\Vert^{3/2}\Vert y\Vert^{3/2}}.
	\]
	The above shows that $M(\beta)\geq\E_{x,y}f^*(x,y)f(x,y)$. We wish to apply the inequality in Proposition~\ref{prop:infinitematrix}, so we will need an upper bound on $\rank(f)$. Define $b\colon\Omega^2\to\bh$ by
	\[
	b(x,y):={|\langle x,y\rangle|^2\over \Vert x\Vert^{3/2}\Vert y\Vert^{3/2}}.
	\]
	We first argue that $\rank(f)\leq\rank(b)$.
	
	Set $r=\rank(b)$ (it is clear that $r<\infty$), and let $g_i,h_i\colon \Omega\to\bh$, $i\in[r]$, be functions so that $b(x,y)=\sum_{i=1}^r g_i(x)h_i(y)$. Now, define functions $s_i,t_i$ by
	\begin{align*}
		s_i(x) &:= g_i(x)+\beta\cdot k^{1/2}\cdot\Vert x\Vert^{1/2}\cdot\E_z\bigl( \Vert z\Vert^{3/2}g_i(z)\bigr),\text{ and}\\
		t_i(y) &:= h_i(y)-\beta\cdot k^{1/2}\cdot\Vert y\Vert^{1/2}\cdot\E_z\bigl( \Vert z\Vert^{3/2}h_i(z)\bigr).	\end{align*}
	We start by noting that for any fixed $x,y$,
	\begin{align*}
		\sum_{i=1}^r g_i(x)\biggl(\Vert y\Vert^{1/2}\cdot\E_z \bigl(\Vert z\Vert^{3/2}h_i(z)\bigr)\biggr) &=\Vert y\Vert^{1/2}\E_z\biggl(\Vert z\Vert^{3/2}\sum_{i=1}^r g_i(x)h_i(z)\biggr)\\
		&= {\Vert y\Vert^{1/2}\over \Vert x\Vert^{3/2}}\E_z|\langle x,z\rangle|^2={1\over k}\Vert x\Vert^{1/2}\Vert y\Vert^{1/2},
	\end{align*}
	as $\mu$ is isotropic. Similarly,
	\[
	\sum_{i=1}^r h_i(y)\biggl(\Vert x\Vert^{1/2}\E_z\bigl(\Vert z\Vert^{3/2}g_i(z)\bigr)\biggr)={1\over k}\Vert x\Vert^{1/2}\Vert y\Vert^{1/2}.
	\]
	Using this, we calculate,
	\begin{align*}
		\sum_{i=1}^r s_i(x)t_i(y) &= \sum_{i=1}^r g_i(x)h_i(y)-\beta^2 k\Vert x\Vert^{1/2}\Vert y\Vert^{1/2}\E_{z,w}\biggl(\Vert z\Vert^{3/2}\Vert w\Vert^{3/2}\sum_{i=1}^r g_i(z)h_i(w)\biggr)\\
		&= b(x,y)-\beta^2k\Vert x\Vert^{1/2}\Vert y\Vert^{1/2}\E_{z,w}|\langle z,w\rangle|^2\\
		&= b(x,y)-\beta^2\Vert x\Vert^{1/2}\Vert y\Vert^{1/2}\\
		&=(1+\beta)f(x,y).
	\end{align*}
	Hence, $\rank(f)\leq r=\rank(b)$, so we only need an upper bound on $\rank(b)$. Here, we break into cases depending on whether $\bh=\R$ or $\bh=\C$. Define $c\colon \Omega^2\to\bh$ by
	\[
	c(x,y):={\langle x,y\rangle\over\Vert x\Vert^{3/4}\Vert y\Vert^{3/4}},
	\]
	which has $\rank(c)=k$.	
		
	\emph{Case 1.} $\bh=\R$. In this case, $b=c^2$, so by Proposition~\ref{prop:hadamard}, we have $\rank(b)\leq{k+1\choose 2}$, which gives the same inequality on $\rank(f)$. Thus, applying Proposition~\ref{prop:infinitematrix}, we bound 
	\begin{align*}
	M(\beta) &\geq \biggl(\E_x{|\langle x,x\rangle|^2-\beta^2\Vert x\Vert^2\Vert x\Vert^2\over (1+\beta)\Vert x\Vert^{3/2}\Vert x\Vert^{3/2}}\biggr)^2\big/{k+1\choose 2}\\
	&=\biggl(\E_x{\Vert x\Vert(1-\beta^2)\over 1+\beta}\biggr)^2\big/{k+1\choose 2}\\
	&=(1-\beta)^2\bigl(\E_x\Vert x\Vert\bigr)^2\big/{k+1\choose 2}.
	\end{align*}
	Combining this lower bound on $M(\beta)$ with the upper bound in Equation~\eqref{eqn:upper}, we have
	\[
	2\beta\alpha\leq{1\over k}+\biggl(\beta^2-{(1-\beta)^2\over{k+1\choose 2}}\biggr)\bigl(\E_x\Vert x\Vert\bigr)^2,
	\]
	for all $\beta\geq 0$. Selecting $\beta=1/\sqrt{k+2}$, we calculate
	\begin{align*}
	{2\alpha\over\sqrt{k+2}} &\leq {1\over k}+\biggl({1\over k+2}-{2(\sqrt{k+2}-1)^2\over k(k+1)(k+2)}\biggr)\bigl(\E_x\Vert x\Vert\bigr)^2\\
	&\leq{1\over k}+{1\over k+2}-{2(\sqrt{k+2}-1)^2\over k(k+1)(k+2)},
	\end{align*}
	where the last line holds because ${1\over k+2}\geq{2(\sqrt{k+2}-1)^2\over k(k+1)(k+2)}$ for all $k\geq 1$ and $\bigl(\E_x\Vert x\Vert\bigr)^2\leq\E_x\Vert x\Vert^2=1$. Solving for $\alpha$ in this expression yields
	\[
	\E_{x,y}|\langle x,y\rangle|=\alpha\leq{(k-1)\sqrt{k+2}+2\over k(k+1)}.
	\]
	
	\emph{Case 2.} $\bh=\C$. Here we have $b=\overline{c}c$, so by Proposition~\ref{prop:hadamard}, we know that $\rank(b)\leq k^2$. Applying Proposition~\ref{prop:infinitematrix} and following the same steps as in Case 1 shows
	\[
	M(\beta)\geq (1-\beta)^2\bigl(\E_x\Vert x\Vert\bigr)^2/k^2\implies 2\beta\alpha\leq{1\over k}+\biggl(\beta^2-{(1-\beta)^2\over k^2}\biggr)\bigl(\E_x\Vert x\Vert\bigr)^2.
	\]
	In this case, we select $\beta=1/\sqrt{k+1}$, which yields
	\begin{align*}
	{2\alpha\over\sqrt{k+1}} &\leq{1\over k}+\biggl({1\over k+1}-{(\sqrt{k+1}-1)^2\over k(k+1)}\biggr)\bigl(\E_x\Vert x\Vert\bigr)^2\\
	&\leq {1\over k}+{1\over k+1}-{(\sqrt{k+1}-1)^2\over k(k+1)},
	\end{align*}
	and solving for $\alpha$ gives
	\[
	\E_{x,y}|\langle x,y\rangle|=\alpha\leq{(k-1)\sqrt{k+1}+1\over k^2}.
	\]\medskip
	
	We now look at the case of equality. 
	
	Let $\alpha(\R)={(k-1)\sqrt{k+2}+2\over k(k+1)}$, $\beta(\R)=1/\sqrt{k+2}$ and $N(\R)={k+1\choose 2}$. Also let $\alpha(\C)={(k-1)\sqrt{k+1}+1\over k^2}$, $\beta(\C)=1/\sqrt{k+1}$ and $N(\C)=k^2$. The proof is identical over $\R$ and $\C$ except for the values of these parameters, so for $\bh\in\{\R,\C\}$, set $\alpha=\alpha(\bh)$, $\beta=\beta(\bh)$ and $N=N(\bh)$. Notice that $\alpha=\beta+(1-\beta)/N$.
	
	First, we establish the ``if'' direction. Let $X$ be a system of $N$ equiangular lines in $\bh^k$. 
It is known\footnote{See, for instance, \cite{lemmens1973equiangular}. We will also re-derive this in the proof of Theorem~\ref{thm:equiconst}; see Equation~\eqref{eqn:angle}.}  that 
for any $x\neq y\in X$, $|\langle x,y\rangle|=\beta$.
Will show in the proof of Theorem~\ref{thm:equiconst}, in Equation~\eqref{eqn:equiso}, that any probability mass $\mu$ on $\bh^k$ with $\mu(x)+\mu(-x)=1/N$ for all $x\in X$ is indeed isotropic. Fix such a mass~$\mu$.
	We calculate,
	\[
	\E_{x,y}|\langle x,y\rangle|=\beta+(1-\beta)\Pr[x\in\{\pm y\}]=\beta+{1-\beta\over N}=\alpha.
	\]

	Now, for the ``only if'' direction, suppose that $\mu$ is isotropic and $\E_{x,y}|\langle x,y\rangle|=\alpha$. Thus, every inequality in the proof of the upper bound must hold with equality. From these equalities, we know the following:
	
	\begin{itemize}
		\item $\Pr[x=0]=0$, otherwise we could construct an isotropic probability mass $\mu'$ with $\E_{x,y\sim\mu'}|\langle x,y\rangle|>\E_{x,y\sim\mu}|\langle x,y\rangle|$, as we showed at the beginning of the proof.
		
		\item\label{itm:unit} If equality holds in Equation~\eqref{eqn:norm}, then it must be the case that $\Vert x\Vert=\Vert y\Vert$ for $\mu$-a.e.\ $x,y\in\bh^k$. As $\mu$ is isotropic, we have $\E_x\Vert x\Vert^2=1$, so we know that $\Vert x\Vert=1$ for $\mu$-a.e.\ $x\in\bh^k$.
		
		\item\label{itm:equiang} If equality holds in Equation~\eqref{eqn:CS}, then it must be the case that for $\mu$-a.e.\ $x,y\in\bh^k$, we have $|\langle x,y\rangle|\in\{\Vert x\Vert\Vert y\Vert,\beta\Vert x\Vert\Vert y\Vert\}$. Since $\Vert x\Vert=1$ for $\mu$-a.e.\ $x\in\bh^k$, it follows that for $\mu$-a.e.\ $x,y\in\bh^k$,
		\[
		|\langle x,y\rangle|=\begin{cases}
		1 & \text{if $x\in\{\pm y\}$},\\
		\beta & \text{otherwise}.
		\end{cases}
		\]
	\end{itemize}

	Therefore, $\supp(\mu)\subseteq X\cup(-X)$ where $X\subseteq\bh^k$ is a system of equiangular lines with $|\langle x,y\rangle|=\beta$ for all $x\neq y\in X$; in particular, $|X|\leq N$.
	
	Recalling that $\alpha=\beta+{1-\beta\over N}$,
	\[
	\beta+{1-\beta\over N}=\E_{x,y}|\langle x,y\rangle|=\beta+(1-\beta)\Pr[x\in\{\pm y\}]\geq\beta+{1-\beta\over|X|}.
	\]
	Therefore, $|X|\geq N$ as well, so $X$ is a system of $N$ equiangular lines over $\bh^k$. Additionally, as $|X|=N$, this means that the inequality above is in fact equality, so $\mu(x)+\mu(-x)=1/N$ for every $x\in X$, as claimed.
\end{proof}

\paragraph{Putting everything together.} We are now ready to give upper bounds on $\SLone_\bh(d,k)$ and analyze the case of equality. To do this, it will be important to know that $\SLone_\bh(d,k)$ is actually achieved.

\begin{prop}\label{prop:supachieved}
	For $\bh\in\{\R,\C\}$ and all positive integers $d,k$, there is some $\mu\in\P_\bh(d,k)$ with $\Lone_\bh(\mu)=\SLone_\bh(d,k)$. 
	
\end{prop}
\begin{proof}
	Let $\{\mu_n\in\P_\bh(d,k):n\in\Z^+\}$ be such that $\SLone_\bh(d,k)\leq\Lone_\bh(X_n)+1/n$ for every $n\in\Z^+$. By Proposition~\ref{prop:iso}, we may suppose that $\mu_n$ is isotropic for all $n\in\Z^+$. As $\mu_n\in\P_\bh(d,k)$, let $X_n=\{x_1^n,\dots,x_{d+k}^n\}$ be a (multi)set so that $\mu_n$ is the uniform distribution over $X_n$. Since $\mu_n$ is isotropic, we know that ${1\over d+k}\sum_{i=1}^{d+k}\Vert x_i^n\Vert^2=\E_{x\sim\mu_n}\Vert x\Vert^2=1$, so it must be the case that $\Vert x_i^n\Vert^2\leq d+k$ for every $i\in[d+k]$ and $n\in\Z^+$. As such, for each $i\in[d+k]$, the sequence $\{x_i^n\}_{n=1}^\infty$ is bounded, so it has a convergent subsequence. Hence, without loss of generality, we may suppose that $\{x_i^n\}_{n=1}^\infty$ converges for every $i\in[d+k]$ and set $x_i=\lim_{n\to\infty}x_i^n$. Let $X=\{x_1,\dots,x_{d+k}\}$ and let $\mu$ be the uniform distribution over $X$. We claim that $\mu$ is isotropic. Indeed, as each $\mu_n$ is isotropic, for any $v\in\bh^k$, we have 
	\[
	\E_{x\sim\mu}|\langle v,x\rangle|^2=\lim_{n\to\infty}\E_{x\sim\mu_n}|\langle v,x\rangle|^2={1\over k}\Vert v\Vert^2.
	\]
	As $\mu$ is isotropic, it must be the case that $\supp(\mu)$ spans $\bh^k$, so as $X$ is a (multi)set of $d+k$ vectors, we have $\mu\in\P_\bh(d,k)$. Now, fix any $i\in[d+k]$ so that $x_i\neq 0$ and any $v\in\bh^k\setminus\{0\}$. We find
	\[
	{\E_{x\sim\mu}|\langle v,x\rangle|\over|\langle v,x_i\rangle|}=\lim_{n\to\infty}{\E_{x\sim\mu_n}|\langle v,x\rangle|\over|\langle v,x_i^n\rangle|}\geq\lim_{n\to\infty}\Lone_\bh(\mu_n)\geq\lim_{n\to\infty}\biggl(\SLone_\bh(d,k)-{1\over n}\biggr)=\SLone_\bh(d,k).
	\]
	Thus $\Lone_\bh(\mu)=\SLone_\bh(d,k)$.
\end{proof}

With this out of the way, we are ready to bound $\SLone_\bh(d,k)$.

\begin{theorem}\label{thm:isobound}\hspace{2em}\vspace{-1ex}
	\begin{enumerate}[label=(\alph*)]
		\item For positive integers $d,k$, 
		\[
		\SLone_\R(d,k)\leq{(k-1)\sqrt{k+2}+2\over k(k+1)},
		\]
		and if equality holds, then there exist ${k+1\choose 2}$ equiangular lines in $\R^k$ and $d\equiv -k\pmod{{k+1\choose 2}}$.
		
		\item For positive integers $d,k$,
		\[
		\SLone_\C(d,k)\leq{(k-1)\sqrt{k+1}+1\over k^2},
		\]
		and if equality holds, then there exist $k^2$ equiangular lines in $\C^k$ and $d\equiv -k\pmod{k^2}$.
	\end{enumerate}

\end{theorem}

In Section~\ref{app:alt}, we give a very different proof that $\SLone_\R(d,2)\leq{2\over 3}$, which may be of separate interest. This alternative proof works by circumscribing an affine copy of a regular hexagon and does not use Lemma~\ref{lem:iso}.

\begin{proof}
	Let $\bh\in\{\R,\C\}$ and suppose $\SLone_\bh(d,k)=\alpha$. By Proposition~\ref{prop:supachieved}, we can find $\mu\in\P_\bh(d,k)$ with $\Lone_\bh(\mu)=\alpha$; we may suppose $\mu$ is isotropic by Proposition~\ref{prop:iso}. As $\Lone_\bh(\mu)=\alpha$, for every $v\in\bh^k$ and $y\in\supp(\mu)$, we must have $\E_x|\langle x,v\rangle|\geq\alpha|\langle y,v\rangle|$. By selecting $v=y$ and averaging over all $y\in\supp(\mu)$, this implies that
	\[
	\E_{x,y}|\langle x,y\rangle|\geq \alpha\E_y|\langle y,y\rangle|=\alpha\E_y\Vert y\Vert^2=\alpha,
	\]
	where the last equality follows from the fact that $\mu$ is isotropic. Lemma~\ref{lem:iso} then gives the upper bound on $\alpha=\SLone_{\bh}(d,k)$.
	
	If $\bh=\R$ and equality holds, then as $\mu$ is isotropic, by Lemma~\ref{lem:iso}, there is a system of ${k+1\choose 2}$ equiangular lines $X\subseteq\R^k$ so that $\mu(x)+\mu(-x)=1/{k+1\choose 2}$ for every $x\in X$, in particular, such a system of equiangular lines must exist. Since $\mu\in\P_\R(d,k)$, we know that $(d+k)\mu(x)\in\Z$ for all $x\in\R^k$, so we must have $(d+k)/{k+1\choose 2}\in\Z$, so $d\equiv -k\pmod{{k+1\choose 2}}$.
	
	The claim is established similarly when $\bh=\C$.
\end{proof}

Theorem~\ref{thm:genlower} follows by combining Theorems~\ref{thm:probbound} and~\ref{thm:isobound}.

\section{Upper bounds}\label{sec:upper}

In this section, we present constructions that yield upper bounds on $\theta_\bh(d,k)$.

We start by proving a general theorem which shows that in order to upper bound $\inner_\bh(d,k)$ it suffices to find an appropriate matrix. For a Hermitian matrix~$C$, denote the largest eigenvalue of $C$ by $\lambda_{\max}(C)$.

\begin{lemma}\label{lem:genconst}
	For $\bh\in\{\R,\C\}$, Let $C\in\bh^{n\times n}$ be Hermitian with $C_{ii}=1$ and $|C_{ij}|\leq 1$ for all $i,j$. If~$\lambda_{\max}(C)$ has multiplicity~$k$ and $d\equiv-k\pmod n$, then
	\[
	\inner_\bh(d,k)\leq{n\over \lambda_{\max}(C)\cdot(d+k)-n}.
	\]
\end{lemma}

\begin{proof}
	As $d\equiv -k\pmod n$, let $b$ be so that $d=nb-k$. Set $\lambda=\lambda_{\max}(C)$ and set $\epsilon={1\over b\lambda-1}$, so $1+\epsilon=\epsilon b\lambda$. It is important to note that $\epsilon>0$. Indeed, if $C\neq I_n$, then as $\tr(C)=n$, we must have $\lambda_{\max}(C)>1$. If it happens to be the case that $C=I_n$, then $k=n$, so as $d>0$, we have $b\geq 2$.
	
	Consider the matrix $A:=(1+\epsilon) I_{nb}-\epsilon(C\otimes J_b)$, where $\otimes$ is the Kronecker/tensor product and $J_b$ is the $b\times b$ all-ones matrix. Note that $A$ is Hermitian, and $A\in\bh^{(d+k)\times (d+k)}$.
	
	As $\lambda=\lambda_{\max}(C)$ has multiplicity $k$, let $N\in\bh^{n\times k}$ have $\rank(N)=k$ and $CN=\lambda N$. Thus $N\otimes J_b$ also has rank $k$ and
	\[
	A(N\otimes J_b)=(1+\epsilon)(N\otimes J_b)-\epsilon(C\otimes J_b)(N\otimes J_b)=(1+\epsilon-\epsilon b\lambda)(N\otimes J_b)=0,
	\]
	by the choice of $\epsilon$. As such, $\rank(A)\leq nb-k=d$. Furthermore, as $\lambda=\lambda_{\max}(C)$, we observe that $A$ is positive semidefinite. Additionally, as $C_{ii}=1$ and $|C_{ij}|\leq 1$ for all $i,j$, we have $A_{ii}=1$ and $|A_{ij}|\leq\epsilon$ for all $i\neq j$.  Therefore,
	\[
	\inner_\bh(d,k)\leq\off(A)\leq\epsilon={1\over b\lambda-1}={n\over\lambda\cdot(d+k)-n}.\qedhere
	\]
\end{proof}

Motivated by the reduction to isotropic measures in the previous section, our usage of Lemma~\ref{lem:genconst} will roughly go as follows: we look for unit vectors $x_1,\dots,x_n\in\bh^k$ so that $|\langle x_i,x_j\rangle|$ is small for all $i\neq j$ and the vectors are, up to scaling, in isotropic position; that is to say $\sum_i x_ix_i^*=\lambda I_k$ for some $\lambda\in\R^+$. In this case, if $A=[x_1|\cdots|x_n]\in\bh^{k\times n}$, we know that $A^*A$ has $1$'s on the diagonal and small entries off the diagonal. Furthermore, $AA^*=\sum_i x_ix_i^*=\lambda I_k$, so $A^*A$ has eigenvalues $\lambda$, with multiplicity $k$, and $0$, with multiplicity $n-k$. We will then let $C$ be an appropriately scaled version of $A^*A$ and apply Lemma~\ref{lem:genconst}.

At this point, it is pertinent to mention that collections of vectors which satisfy $\sum_i x_ix_i^*=\lambda I_k$ for some $\lambda\in\R^+$ are also known as \emph{finite tight frames} and have been studied extensively in the literature (see~\cite{Christensen2016,kovavcevic2008introduction} for a survey). We will rely on known constructions of finite tight frames.

We will be able to execute the above plan for only some values of $d$ and $k$; we deal with the remaining values using monotonicity of $\inner$, that is
$\inner_\bh(d,k)\leq\inner_\bh(d',k')$ whenever $d\geq d'$ and $k\leq k'$.

We can now apply this general construction to prove Theorem~\ref{thm:equiconst}; namely, if large systems of equiangular lines exist, then the lower bound in Theorem~\ref{thm:genlower} is tight.

\begin{proof}[Proof of Theorem~\ref{thm:equiconst}]
	Theorem~\ref{thm:genlower} establishes the lower bound for all $d,k$, so we need establish only the upper bound.
	
	Let $\{x_1,\dots,x_N\}\subseteq\bh^k$ be a system of equiangular lines where $N={k+1\choose 2}$ (if~$\bh=\R$) or $N=k^2$ (if~$\bh=\C$). From Gerzon's proof\footnote{Gerzon never published his proof. The original reference appears to be \cite[Theorem~3.5]{lemmens1973equiangular}. The proof can also be found in \cite[Miniature~9]{matousek2010thirty}.} that there are at most $N$ equiangular lines in $\bh^k$, we know that the projection matrices $x_1x_1^*,\dots,x_Nx_N^*$ span the space of all Hermitian matrices in $\bh^{k\times k}$ as a vector space over $\R$. Thus, there are constants $c_1,\dots,c_N\in\R$ for which $I_k=\sum_i c_ix_ix_i^*$. Let $\beta$ be the common inner product of $\{x_1,\dots,x_N\}$, that is, $|\langle x_i,x_j\rangle|=\beta$ for all $i\neq j$. For any fixed $j\in[N]$,
	\begin{align*}
	1 &=\tr( x_jx_j^*)=\tr( I_kx_jx_j^*)\\
	&=\tr\biggl(\sum_ic_ix_ix_i^*x_jx_j^*\biggr)=\sum_ic_i|\langle x_i,x_j\rangle|^2\\
	&=c_j+\sum_{i\neq j}c_i\beta^2=(1-\beta^2)c_j+\sum_ic_i\beta^2,
	\end{align*}
	so for all $j$,
	\[
	c_j={1-\beta^2\sum_i c_i\over 1-\beta^2}.
	\]
	In particular, $c_1=\dots=c_N=c$. Now,
	\[
	k=\tr (I_k)=\tr\biggl(\sum_i c x_ix_i^*\biggr)=cN,
	\]
	so $c={k\over N}$. Hence,
	\begin{equation}\label{eqn:equiso}
	\sum_i x_ix_i^*={N\over k}I_k,
	\end{equation}
	 and
	\begin{equation}\label{eqn:angle}
	{k\over N}={1-\beta^2\sum_{i=1}^N{k\over N}\over 1-\beta^2}\implies \beta=\sqrt{{N-k\over kN-k}}.
	\end{equation}
	
	Now, let $A=[x_1|\cdots|x_N]$, so $(A^*A)_{ii}=1$ and $|(A^*A)_{ij}|=\beta$ for all $i\neq j$. Additionally, $AA^*=\nobreak \sum_ix_ix_i^*={N\over k}I_k$, so $A^*A$ has eigenvalues ${N\over k}$ and $0$, where the former has multiplicity~$k$. Finally, set $C:={1\over\beta}(A^*A-I_N)+I_N\in\bh^{N\times N}$, which is Hermitian with $C_{ii}=1$ and $|C_{ij}|=1$ for all $i,j$. Furthermore, $\lambda_{\max}(C)={1\over\beta}\bigl({N\over k}-1\bigr)+1$, which has multiplicity~$k$. 
	
	If $\bh=\R$, substituting $N={k+1\choose 2}$ shows that ${\lambda_{\max}(C)\over N}=\alpha_k$. By Lemma~\ref{lem:genconst}, if $d\equiv -k\pmod{{k+1\choose 2}}$, 
	\[
	\inner_\R(d,k)\leq{N\over\lambda_{\max}(C)(d+k)-N}={1\over\alpha_k(d+k)-1}.
	\]
	
	If $\bh=\C$, substituting $N=k^2$ shows that ${\lambda_{\max}(C)\over N}=\alpha_k^*$. By Lemma~\ref{lem:genconst}, if $d\equiv -k\pmod{k^2}$,
	\[
	\inner_\C(d,k)\leq{N\over\lambda_{\max}(C)(d+k)-N}={1\over\alpha_k^*(d+k)-1}.\qedhere
	\]
\end{proof}

For $k\in\{1,2,3,7,23\}$, there are in fact systems of ${k+1\choose 2}$ equiangular lines over $\R^k$, so in these cases, we can pin down $\inner_\R(d,k)$ precisely for infinitely many values of~$d$. We have previously mentioned the value of $\inner_\R(d,1)$ in Corollary~\ref{cor:k=1}, so we do not restate it here.
\begin{cor}\label{cor:equicases}\hspace{2em}\vspace{-1ex}
	\begin{itemize}
		\item If $d\equiv -2\pmod{3}$, then $\off_\R(d,2)=\inner_\R(d,2)={3\over 2d+1}$.
		\item If $d\equiv -3\pmod{6}$, then $\off_\R(d,3)=\inner_\R(d,3)={6\over (\sqrt{5}+1)d+3(\sqrt{5}-1)}$.
		\item If $d\equiv -7\pmod{28}$, then $\off_\R(d,7)=\inner_\R(d,7)={14\over 5d+21}$.
		\item If $d\equiv -23\pmod{276}$, then $\off_\R(d,23)=\inner_\R(d,23)={69\over 14d+253}$.
	\end{itemize}
	
\end{cor}

Over $\C$, the existence of $k^2$ equiangular lines over $\C^k$ is known for numerous values of $k$. For example, constructions exist for $k\in\{1,2,\dots,16,19,24,28,35,48\}$, and, up to numerical precision, all $k\leq 67$ (see~\cite{Wiebe_2013} for a survey). In fact, it is conjectured that there are $k^2$ equiangular lines over $\C^k$ \emph{for all} $k$. Thus, conjecturally, we have the following:
\begin{conj}
	For every positive integer $k$, if $d\equiv -k\pmod{k^2}$, then
	\[
	\inner_\C(d,k)={1\over \alpha_k^*(d+k)-1},
	\]
	where $\alpha_k^*={(k-1)\sqrt{k+1}+1\over k^2}$.
\end{conj}

We now turn to upper bounds on $\inner_\bh(d,k)$ in the case when no system of equiangular lines of size $\binom{k+1}{2}$ (if $\bh=\R$) or $k^2$ (if $\bh=\C$) exists.
\begin{defn}
	For $\bh\in\{\R,\C\}$, matrices $B_1,\dots,B_\ell\in\bh^{k\times k}$ are said to be \emph{mutually unbiased bases of $\bh^k$} if $B_i^*B_i=I_k$ for all $i$ and every entry $B_i^*B_j$ has magnitude $1/\sqrt{k}$ for all $i\neq j$. 
\end{defn}
	
	The following is known:
	\begin{itemize}
		\item If $k$ is a power of $4$, then there is a collection of ${k\over 2}+1$ mutually unbiased bases of $\R^k$ (see~\cite{cameron1973quadratic}). 
		\item If $k$ is a prime power, then there is a collection of $k+1$ mutually unbiased bases of $\C^k$ (see~\cite{bengtsson2007three}). 
	\end{itemize}

\begin{lemma}\label{lem:mubconst}
	For $\bh\in\{\R,\C\}$, if there exists a collection of $\ell$ mutually unbiased bases of $\bh^k$, then whenever $d\equiv -k\pmod{k\ell}$, 
	\[
	\inner_\bh(d,k)\leq{k\ell\over\bigl(\sqrt{k}(\ell-1)+1\bigr)(d+k)-k\ell}.
	\]
\end{lemma}
\begin{proof}
	Let $B_1,\dots,B_\ell$ be a collection of mutually unbiased bases of $\bh^k$ and consider the matrix $A=[B_1|B_2|\cdots|B_\ell]$. From the properties of mutually orthogonal bases, we find that $AA^*=\ell I_k$, so $A^*A$ has eigenvalues $\ell$ and $0$ where the former has multiplicity $k$. Furthermore, $A^*A$ has $1$'s on the diagonal and every off-diagonal entry is either $0$ or has magnitude $1/\sqrt{k}$. Set $C=\sqrt{k}(A^*A-I_{k\ell})+I_{k\ell}$, so $C\in\bh^{k\ell\times k\ell}$ is a Hermitian matrix with $C_{ii}=1$ and $|C_{ij}|\leq 1$ for all $i,j$. Additionally, $\lambda_{\max}(C)=\sqrt{k}(\ell-1)+1$, which has multiplicity $k$, so the claim follows by applying Lemma~\ref{lem:genconst}.
\end{proof}

Using the above lemma, we can prove Theorem~\ref{thm:genupper} over $\R$ for infinitely many values of $k$ and give a bound that is off by a factor of at most~$2$ for general~$k$.

\begin{theorem}\label{thm:offby2}
	If $k$ is a power of $4$, then whenever $d\equiv -k\pmod{k^2/2+k}$,
	\[
	\inner_\R(d,k)\leq{\sqrt{k+4-\Omega(k^{-1/2})}\over d}.
	\]
	Additionally, for any fixed $k$,
	\[
	\inner_\R(d,k)\leq \bigl(1+o(1)\bigr){2\sqrt{k+1}\over d}.
	\]
\end{theorem}
\begin{proof}
	If $k$ is a power of $4$, then there is a collection of $\ell={k\over 2}+1$ mutually unbiased bases of $\R^k$. Thus, by Lemma~\ref{lem:mubconst}, whenever $d\equiv -k\pmod{k^2/2+k}$,
	\[
	\inner_\R(d,k)\leq{k^2/2+k\over\bigl(k^{3/2}/2+1\bigr)(d+k)-k^2/2-k}\leq{\sqrt{k+4-\Omega(k^{-1/2})}\over d}.
	\]
	
	For a general $k$, let $k'$ be a power of $4$ satisfying $k\leq k'\leq 4k$. By monotonicity,
	\[
	\inner_\R(d,k)\leq \inner_\R(d,k')\leq\bigl(1+o(1)\bigr){2\sqrt{k+1}\over d},
	\]
	as $d\to\infty$.
\end{proof}

In the case of complex numbers, we can establish Theorem~\ref{thm:genupper} immediately.

\begin{theorem}\label{thm:asymtightcomp}
	If $q$ is a prime power, then whenever $d\equiv -q\pmod{q^2+q}$,
	\[
	\inner_\C(d,q)\leq {\sqrt{q+2-\Omega(q^{-1/2})}\over d}.
	\]
	Additionally, for any fixed $k$,
	\[
	\inner_\C(d,k)\leq\bigl(1+o(1)\bigr){\sqrt{k+O(k^{21/40})}\over d}.
	\]
\end{theorem}
\begin{proof}
	If $q$ is a prime power, then there is a collection of $\ell=q+1$ mutually unbiased bases of $\C^k$. Thus, whenever $d\equiv -q\pmod{q^2+q}$,
	\[
	\inner_\C(d,q)\leq{q^2+q\over \bigl(q^{3/2}+1\bigr)(d+q)-q^2-q}\leq{\sqrt{q+2-\Omega(q^{-1/2})}\over d}
	\]
	
	For any $k$, since there is always some prime $q$ satisfying $k\leq q\leq k+O(k^{21/40})$ (see~\cite{baker2001difference}), by monotonicity, we have
	\[
	\inner_\C(d,k)\leq \inner_\C(d,q)\leq\bigl(1+o(1)\bigr){\sqrt{k+O(k^{21/40})}\over d},
	\]
	as $d\to\infty$.
\end{proof}

Notice that Theorem~\ref{thm:asymtightcomp} implies that there is a constant $c$ such that for any $\epsilon>0$, if $k>c\epsilon^{-40/19}$, then
\[
\inner_\C(d,k)\leq\bigl(1+o(1)\bigr){(1+\epsilon)\sqrt{k}\over d},
\]
which establishes Theorem~\ref{thm:genupper} over the complex numbers.

We now present a more general construction of nearly orthogonal vectors which makes use of Steiner systems and Hadamard matrices. This construction will allow us to establish Theorem~\ref{thm:genupper} over the real numbers.

\begin{defn}
	A $(2,\ell,n)$-Steiner system consists of $n$ points and a collection of subsets of these points, called blocks, where each block contains exactly $\ell$ points and any two points are contained in exactly one block together.\footnote{In standard notation, $k$ is used in place of $\ell$ when discussing Steiner systems, but we opt to go against this in order to stay consistent with the notation in this paper.} If $k$ is the number of blocks and $r$ is the degree of any point, it is well-known that $k={n(n-1)\over\ell(\ell-1)}$ and $r={n-1\over \ell-1}$.
\end{defn}

\begin{defn}
	For $\bh\in\{\R,\C\}$, a Hadamard matrix over $\bh$ of order $n$ is a matrix $H\in\bh^{n\times n}$ so that for all $i,j$, $|H_{ij}|=1$ and $H^*H=nI_n$. When $\bh=\C$, Hadamard matrices of order $n$ exists for all $n$. When $\bh=\R$, it is not known for which $n$ Hadamard matrices of order $n$ exist. It is known however that such an $n>2$ must be divisible by $4$.
\end{defn}

The following tight frame was constructed by Fickus, Mixon and Tremain~\cite{fickus2012steiner}. We state their construction in language which will be useful for our purposes.

\begin{theorem}[Fickus, Mixon and Tremain~{\cite[Theorem 1]{fickus2012steiner}}]\label{thm:fickusteiner}
	Let $\bh\in\{\R,\C\}$ and suppose there exists a $(2,\ell,n)$-Steiner system with $k$ blocks and degree $r$. If, in addition, there exists a Hadamard matrix of order $r+1$ over $\bh$, then there is a matrix $B\in\bh^{k\times n(r+1)}$ satisfying:
	\begin{itemize}
		\item $B^*B$ has $r$'s on the diagonal and every off-diagonal entry has magnitude~$1$, and
		\item $BB^*=\ell(r+1)I_k$.
	\end{itemize}
\end{theorem}

From this construction, we can give bounds on $\inner_\bh(d,k)$.

\begin{cor}\label{cor:steinerconst}
		Let $\bh\in\{\R,\C\}$ and suppose there exists a $(2,\ell,n)$-Steiner system with $k$ blocks and degree $r$. If, in addition, there exists a Hadamard matrix of order $r+1$ over $\bh$, then whenever $d\equiv -k\pmod{n(r+1)}$,
	\[
	\inner_\bh(d,k)\leq{n(r+1)\over \bigl(\ell(r+1)-r+1\bigr)(d+k)-n(r+1)}.
	\]
\end{cor}
\begin{proof}
	Let $B$ be the matrix as in Theorem~\ref{thm:fickusteiner} and set $C:=B^*B-(r-1)I_{n(r+1)}$. We notice that $C\in\bh^{n(r+1)\times n(r+1)}$ is Hermitian with $C_{ii}=1$ and $|C_{ij}|\leq 1$ for all $i,j$. Additionally, as $BB^*=\ell(r+1)I_k$, we know that $\lambda_{\max}(C)=\ell(r+1)-(r-1)$, which has multiplicity $k$. Thus, the claim follows from Lemma~\ref{lem:genconst}.
\end{proof}

Using Corollary~\ref{cor:steinerconst}, we can establish Theorem~\ref{thm:genupper} in the case of the reals.

\begin{theorem}\label{thm:realasym}
	For every $\epsilon>0$, there is a $k_0$ so that whenever $k\geq k_0$,
	\[
	\inner_\R(d,k)\leq\bigl(1+o(1)\bigr){(1+\epsilon)\sqrt{k}\over d}.
	\]
\end{theorem}

In order to prove this, we require the following results:
\begin{fact}[Prime number theorem for arithmetic progressions {\cite[Chapters 20,21]{davenport}}]\label{fact:apprimes}
	For integers $a,n$ with $\gcd(a,n)=1$, there is a function $f_{a,n}$ with $f_{a,n}(x)\to 0$ as $x\to\infty$ so that for any positive $x$, there is a prime $p\equiv a\pmod n$ satisfying $x\leq p\leq \bigl(1+f_{a,n}(x)\bigr)x$.
\end{fact}
\begin{fact}[Wilson~\cite{wilsonSteiner}]\label{fact:steiner}
	For any positive integer $\ell$, there is some other integer $N_\ell$ so that if $n\geq N_\ell$ with $(\ell-1)\mid (n-1)$ and $\ell(\ell-1)\mid n(n-1)$, then a $(2,\ell,n)$-Steiner system exists.
\end{fact}
\begin{fact}[Paley~\cite{Paley_1933}]\label{fact:hadamard}
	Let $q$ be a prime power. If $q\equiv 3\pmod 4$, then a real Hadamard matrix of order $q+1$ exists.
\end{fact}

\begin{proof}[Proof of Theorem~\ref{thm:realasym}]
	If we can locate a prime $p\equiv 3\pmod 4$ and a $(2,\ell,n)$-Steiner system with $k={n(n-1)\over\ell(\ell-1)}$ blocks and degree ${n-1\over\ell-1}=p$, then Fact~\ref{fact:hadamard} and
        Corollary~\ref{cor:steinerconst} together imply that whenever $d\equiv -k\pmod{n(r+1)}$, we have
	\begin{align*}
	\inner_\R(d,k) &\leq{n(r+1)\over\bigl(\ell(r+1)-r+1\bigr)(d+k)-n(r+1)}\\
	&\leq{n(r+1)\over\ell(r+1)-r+1}{1\over d}\\
	&={n\bigl({n-1\over\ell-1}+1\bigr)\over n+\ell}{1\over d}\\
	&=\sqrt{k\biggl(1+{(n-\ell)^2(n+\ell-1)\over (n+\ell)^2(n-1)(\ell-1)}\biggr)}\cdot {1\over d}\\
	&\leq\sqrt{k\biggl(1+{1\over \ell-1}\biggr)}\cdot{1\over d},
	\end{align*}
	where the last line follows as $n\geq\ell$.
	
	Given an $\epsilon>0$, pick $\ell$ odd so that $1/(\ell-1)<\epsilon$. Consider any prime $p$ so that $p\equiv 3\pmod 4$ and $p\equiv 1\pmod \ell$ and $(\ell-1)p\geq N_\ell$, where $N_\ell$ is as in Fact~\ref{fact:steiner}. Set $n=1+(\ell-1)p$.
	
	We notice that ${n-1\over\ell-1}=p$, and also that $k':={n(n-1)\over\ell(\ell-1)}=p\bigl(p-{p-1\over\ell}\bigr)$ is an integer by the choice of~$p$. By Fact~\ref{fact:steiner}, there exists a $(2,\ell,n)$-Steiner system with $k'$ blocks and degree $p$. By Fact~\ref{fact:apprimes}, for any sufficiently large $k$, we can find a suitable prime $p$ for which $k\leq k'\leq(1+\epsilon)k$, so by monotonicity and the remark above,
	\[
	\inner_\R(d,k)\leq\inner_\R(d,k')\leq\bigl(1+o(1)\bigr){(1+\epsilon)\sqrt{k}\over d},
	\]
	as $d\to\infty$.	
\end{proof}

\section{An alternative proof that \texorpdfstring{$\SLone_\R(d,2)\leq{2\over 3}$}{SL\textunderscore R(d,2)<=2/3}}\label{app:alt}

Here we present an alternative proof of the upper bound on $\SLone_\R(d,2)$. We have been unable to generalize this proof to get a bound on $\SLone_\R(d,k)$ for any other $k$.

The proof hinges on the following result, which was proved by Go{\l}ab~\cite{golab1932some} and refined by Besicovitch~\cite{Besicovitch1948}.

\begin{lemma}
	If $C\subseteq\R^2$ is compact, convex and centrally-symmetric, and $H$ is a centrally-symmetric regular hexagon, then there is $Q\in\GL_2(\R)$ so that $QH$ circumscribes $C$.
\end{lemma}

\begin{proof}[Proof that $\SLone_\R(d,2)\leq{2\over 3}$]
	Suppose $\mu\in\P_\R(d,2)$, and let $C$ be the convex hull of $\supp(\mu)\cup(-\supp(\mu))$, so as $\supp(\mu)$ is finite, we know that $C$ is compact, convex and centrally-symmetric.  Let $H$ be the hexagon centered at the origin with distance $2$ between its parallel edges, as shown in Figure~\ref{fig:hexagon}. By the lemma, there is $Q\in\GL_2(\R)$ such that $QH$ circumscribes~$C$. We label the top three lines bounding $H$ as $\ell_1,\ell_2,\ell_3$ where $\ell_i=\{x\in\R^2:\langle x,v_i\rangle=1\}$.

	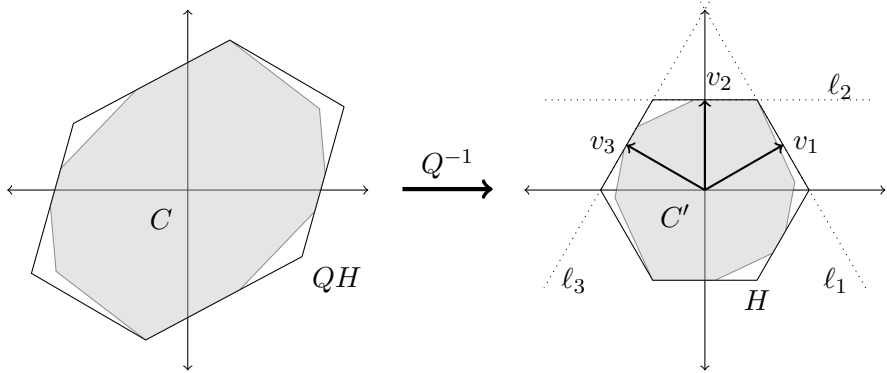
\begin{figure}[ht]
		\centering
		\begin{tikzpicture}[extended line/.style={shorten >=-#1,shorten <=-#1}, extended line/.default=1cm, scale=1.2]
		\def\r{{2/sqrt(3)}}
		\draw[<->] (-2,0)--(2,0);
		\draw[<->](0,-2)--(0,2);
		\pgftransformcm{1.5}{.8}{-0.4}{1.2}{\pgfpoint{0}{0}}
		\foreach \x in {0,...,5}{
			\coordinate (1\x) at (60*\x:\r);}
		\draw (10)--(11)--(12)--(13)--(14)--(15)--(10);
		
		\coordinate (c1) at (60:\r);
		\coordinate (c2) at (5:1);
		\coordinate (c3) at ($(10)!0.45!(15)$);
		\coordinate (c4) at ($(10)!0.7!(15)$);
		\coordinate (c5) at ($(15)!0.4!(14)$);
		\coordinate (c6) at ($(0,0)-(c1)$);
		\coordinate (c7) at ($(0,0)-(c2)$);
		\coordinate (c8) at ($(0,0)-(c3)$);
		\coordinate (c9) at ($(0,0)-(c4)$);
		\coordinate (c10) at ($(0,0)-(c5)$);
		
		\draw[fill=gray!50,opacity=.4] (c1)--(c2)--(c3)--(c4)--(c5)--(c6)--(c7)--(c8)--(c9)--(c10)--cycle;
		
		\node[below right] (H) at (15) {$QH$};
		\node[below left] (C) at (-.05,-.05) {$C$};
		\end{tikzpicture}\quad
		\begin{tikzpicture}[scale=1.2]
		\draw[ultra thick,->] (-.5,0)--(.5,0);
		\node[above] (1) at (0,0) {$Q^{-1}$};
		\node (2) at (0,-1.91) {};
		\end{tikzpicture}\quad
		\begin{tikzpicture}[extended line/.style={shorten >=-#1,shorten <=-#1}, extended line/.default=1cm, scale=1.2]
		\def\r{{2/sqrt(3)}}
		\draw[<->] (-2,0)--(2,0);
		\draw[<->](0,-2)--(0,2);
		\foreach \x in {0,...,5}{
			\coordinate (1\x) at (60*\x:\r);}
		\draw (10)--(11)--(12)--(13)--(14)--(15)--(10);
		\node[right] (1) at (30:1) {$v_1$};
		\node[above right] (2) at (-.1,1) {$v_2$};
		\node[left] (3) at (150:1) {$v_3$};
		\coordinate (21) at (30:1);
		\coordinate (22) at (90:1);
		\coordinate (23) at (150:1);
		
		\coordinate (c1) at (60:\r);
		\coordinate (c2) at (5:1);
		\coordinate (c3) at ($(10)!0.45!(15)$);
		\coordinate (c4) at ($(10)!0.7!(15)$);
		\coordinate (c5) at ($(15)!0.4!(14)$);
		\coordinate (c6) at ($(0,0)-(c1)$);
		\coordinate (c7) at ($(0,0)-(c2)$);
		\coordinate (c8) at ($(0,0)-(c3)$);
		\coordinate (c9) at ($(0,0)-(c4)$);
		\coordinate (c10) at ($(0,0)-(c5)$);
		
		\draw[fill=gray!50,opacity=.4] (c1)--(c2)--(c3)--(c4)--(c5)--(c6)--(c7)--(c8)--(c9)--(c10)--cycle;

		\draw[->,thick] (0,0)--(21);
		\draw[->,thick] (0,0)--(22);
		\draw[->,thick] (0,0)--(23);
		\draw[dotted, extended line=1.5cm] (10) -- (11);
		\draw[dotted, extended line=1.5cm] (11)--(12);
		\draw[dotted, extended line=1.5cm] (12)--(13);
		\node[left] (31) at (1.7,-1) {$\ell_1$};
		\node[above] (32) at (1.5,.9) {$\ell_2$};
		\node[right] (33) at (-1.7,-1) {$\ell_3$};
		\node[below] (H') at (15) {$H$};
		\node[below left] (C') at (-.05,-.05) {$C'$};
		\end{tikzpicture}
		\caption{(Left) $C$ circumscribed by  $QH$. (Right) The result of applying $Q^{-1}$ to $C$ and $QH$. In this image, $v_1,v_2,v_3$ are unit vectors. \label{fig:hexagon}}
	\end{figure}
	
	Set $\mu'=Q^{-1}\mu$ and $C'=Q^{-1}C$, so $C'$ is the convex hull of $\supp(\mu')\cup(-\supp(\mu'))$ and $H$ circumscribes $C'$. Now, consider the maximization problem:
	\[
	\max_{x\in H}\sum_{i=1}^3|\langle x,v_i\rangle|.
	\]
	As $\sum_{i=1}^3|\langle x,v_i\rangle|$ is a convex function and $H$ is also convex, the maximum occurs at a vertex of $H$. Thus, if $\widehat{x}$ denotes such an optimal solution, without loss, $\widehat{x}\in\ell_1\cap\ell_2$, so $\langle \widehat{x},v_1\rangle=\langle \widehat{x},v_2\rangle=1$ and $\langle \widehat{x},v_3\rangle=0$. We conclude that $\sum_{i=1}^3|\langle x,v_i\rangle|\leq 2$ for every $x\in H$.
	
	Therefore, as $\supp(\mu')\subseteq C'\subseteq H$,
	\[
	\sum_{i=1}^3\E_{x\sim\mu'}|\langle x,v_i\rangle|\leq 2\implies \E_{x\sim\mu'}|\langle x,v_i\rangle|\leq{2\over 3},\text{ for some $i\in[3]$.}
	\]
	Without loss, suppose that $\E_{x\sim\mu'}|\langle x,v_1\rangle|\leq {2\over 3}$.
	Finally, as $H$ circumscribes $C'$, for each edge of $H$, there is some vertex of $C'$ lying on this edge. In other words, there is some $y\in\supp(\mu')$ for which $|\langle y,v_1\rangle|=1$, so
	\[
	\Lone_\R(\mu)=\Lone_\R(\mu')\leq{\E_{x\sim\mu'}|\langle x,v_1\rangle|\over|\langle y,v_1\rangle|}\leq{2\over 3}.\qedhere
	\]
\end{proof}

\section{Concluding remarks and open problems}
\begin{itemize}[leftmargin=*]
	\item Because we rely on the existence of general Steiner systems, the dependence of $k_0$ on $\epsilon$ in Theorem~\ref{thm:realasym} is poor. It would be of interest to improve this dependence.
	
	\item When considering upper bounds, we focused on $\inner_\bh(d,k)$ instead of $\off_\bh(d,k)$. For constructions for the latter, one could rephrase Lemma~\ref{lem:genconst} to read: Let $C\in\bh^{n\times n}$ with $C_{ii}=1$ and $|C_{ij}|\leq 1$ for all $i,j$ and let $\lambda$ be any eigenvalue of $C$ with $\lambda\in\bh$. If $\lambda$ has multiplicity $k$ and $d\equiv -k\pmod n$, then
	\[
	\off_\bh(d,k)\leq\bigg|{n\over\lambda\cdot(d+k)-n}\bigg|.
	\]
	This could lead to improved upper bounds on $\off_\bh(d,k)$ which may not hold for $\inner_\bh(d,k)$.
	
	\item Suppose $k$ is such that no system of ${k+1\choose 2}$ equiangular lines exists in $\R^k$; by how much can the lower bound in Theorem~\ref{thm:genlower} be improved?
	
	\item What is $\inner_\R(d,4)$? Theorem~\ref{thm:offby2} shows that $\inner_\R(d,4)\lessapprox{2.4\over d}$ and Theorem~\ref{thm:genlower} shows that $\inner_\R(d,4)\gtrapprox{2.139\over d}$. It would be interesting to close this gap.
	
	\item How small can $\epsilon$ be so that there is some set of $2d+k$ unit vectors $X\subseteq\R^d$ with $\langle x,y\rangle\leq\epsilon$ for all $x\neq y\in X$? Define $\theta'(d,k)$ to be this smallest $\epsilon$. Certainly $\theta'(d,k)\leq\inner_\R(d,\lceil k/2\rceil)\approx{\sqrt{k/2}\over d}$ for a fixed $k$; however, we have be unable to prove matching lower bounds. Using the linear programming method of Delsarte, Goethals and Seidel~\cite{Delsarte1977}, we can show that $\theta'(d,k)\geq(1-o(1)){k\over d^2}$ for a fixed $k$, but it seems unlikely that such an approach will be able to improve this lower bound. Using a different argument, in a forthcoming paper Balla shows that $\theta'(d,k)\geq c_k/d$.
	
	\item For a matrix $A\in\bh^{n\times n}$ and $p>0$, define $\off^p(A):=\bigl(\sum_{i\neq j}|A_{ij}|^p\bigr)^{1/p}$, i.e.\ the $L^p$ norm of the off-diagonal entries of $A$. We then define $\off_\bh^p(d,k):=\min_A\off^p(A)$ where the minimum is taken over all $A\in\bh^{(d+k)\times(d+k)}$ with $A_{ii}=1$ for all $i$ and $\rank(A)= d$. In this context, we can interpret $\off_\bh(d,k)$ as $\off_\bh^\infty(d,k)$. 
	
	For $2\leq p\leq \infty$, by following the arguments in this paper, one can relate the problem of lower-bounding $\off_\bh^p(d,k)$ to finding upper bounds on $\E_{x,y\sim\mu}|\langle x,y\rangle|^q$ where ${1\over p}+{1\over q}=1$ and $\mu$ is an isotropic probability mass on $\bh^k$. We conjecture the following:
	\begin{conj}
		For a positive integer $k$, set $\beta=1/\sqrt{k+2}$ and $N={k+1\choose 2}$. If $1\leq q\leq 2$ and $\mu$ is an isotropic probability mass on $\R^k$, then
		\[
		\E_{x,y\sim\mu}|\langle x,y\rangle|^q\leq \beta^q+{1-\beta^q\over N},
		\]
		with equality if and only if there is $X\subseteq\R^k$, a system of $N$ equiangular lines, and $\mu$ satisfies $\mu(x)+\mu(-x)= 1/N$ for every $x\in X$.
	\end{conj}
	We conjecture also the natural analogue when $\R$ is replaced by $\C$.

\textbf{Added in a revision.} The conjecture has since been resolved by Glazyrin \cite[Theorem~4]{glazyrin_lq}. 
\end{itemize}

\bibliographystyle{abbrv}
\bibliography{references}

\end{document}